\documentclass[a4paper,12pt]{amsart}
\usepackage{amsfonts}
\usepackage{amssymb}
\usepackage{ifthen}
\usepackage{graphicx}
\nonstopmode \numberwithin{equation}{section}
\usepackage[usenames]{color}
\newtheorem{thm}{Theorem}
\newtheorem{lem}{Lemma}
\newtheorem{cor}{Corollary}
\newtheorem{prop}{Proposition}

\newtheorem{examp}{Example}

\newtheorem{cl}{Claim}
\newtheorem{ca}{Case}
\newtheorem{sca}{Subcase}
\newtheorem{scl}{Subclaim}
\newtheorem{conj}[equation]{Conjecture}

\theoremstyle{definition}
\newtheorem{defn}{Definition}
\newtheorem{op}[equation]{Open Problem}
\newtheorem{ques}[equation]{Question}
\newtheorem{rem}{Remark}

\newcounter {own}
\def\theown {\thesection       .\arabic{own}}

\newenvironment{pf}[1][]{%
 \vskip 3mm
 \noindent
 \ifthenelse{\equal{#1}{}}%
  {{\slshape Proof. }}%
  {{\slshape #1.} }%
 }%
{\qed\bigskip}

\newcounter{alphabet}
\newcounter{tmp}
\newenvironment{Thm}[1][]{\refstepcounter{alphabet}%
\bigskip%
\noindent%
{\bf Theorem \Alph{alphabet}}%
\ifthenelse{\equal{#1}{}}{}{ (#1)}%
{\bf .} \itshape}{\vskip 8pt}

\makeatletter
\newcommand{\Ref}[1]{\@ifundefined{r@#1}{}{\setcounter{tmp}{\ref{#1}}\Alph{tmp}}}
\makeatother

\def\be{\begin{equation}}
\def\ee{\end{equation}}

\newcommand{\bee}{\begin{enumerate}}
\newcommand{\eee}{\end{enumerate}}

\newcommand{\blem}{\begin{lem}}
\newcommand{\elem}{\end{lem}}
\newcommand{\bthm}{\begin{thm}}
\newcommand{\ethm}{\end{thm}}
\newcommand{\bcor}{\begin{cor}}
\newcommand{\ecor}{\end{cor}}
\newcommand{\beg}{\begin{examp}}
\newcommand{\eeg}{\end{examp}}
\newcommand{\begs}{\begin{examples}}
\newcommand{\eegs}{\end{examples}}
\newcommand{\bdefe}{\begin{defn}}
\newcommand{\edefe}{\end{defn}}
\newcommand{\bprob}{\begin{prob}}
\newcommand{\eprob}{\end{prob}}
\newcommand{\bques}{\begin{ques}}
\newcommand{\eques}{\end{ques}}
\newcommand{\bei}{\begin{itemize}}
\newcommand{\eei}{\end{itemize}}
\newcommand{\bcon}{\begin{conj}}
\newcommand{\econ}{\end{conj}}
\newcommand{\bop}{\begin{op}}
\newcommand{\eop}{\end{op}}

\newcommand{\bca}{\begin{ca}}
\newcommand{\eca}{\end{ca}}
\newcommand{\bsca}{\begin{sca}}
\newcommand{\esca}{\end{sca}}

\newcommand{\bcl}{\begin{cl}}
\newcommand{\ecl}{\end{cl}}

\newcommand{\bscl}{\begin{scl}}
\newcommand{\escl}{\end{scl}}

\newcommand{\bcons}{\begin{conjs}}
\newcommand{\econs}{\end{conjs}}
\newcommand{\bprop}{\begin{propo}}
\newcommand{\eprop}{\end{propo}}
\newcommand{\br}{\begin{rem}}
\newcommand{\er}{\end{rem}}
\newcommand{\brs}{\begin{rems}}
\newcommand{\ers}{\end{rems}}
\newcommand{\bo}{\begin{obser}}
\newcommand{\eo}{\end{obser}}
\newcommand{\bos}{\begin{obsers}}
\newcommand{\eos}{\end{obsers}}
\newcommand{\bpf}{\begin{pf}}
\newcommand{\epf}{\end{pf}}
\newcommand{\ba}{\begin{array}}
\newcommand{\ea}{\end{array}}
\newcommand{\beq}{\begin{eqnarray}}
\newcommand{\beqq}{\begin{eqnarray*}}
\newcommand{\eeq}{\end{eqnarray}}
\newcommand{\eeqq}{\end{eqnarray*}}

\newcounter{minutes}
\divide\time by 60
\newcounter{hours}
\multiply\time by 60 \addtocounter{minutes}{-\time}

\begin{document}
\bibliographystyle{amsplain}
\title [] {On polyharmonic univalent mappings}

\def\thefootnote{}
\footnotetext{ \texttt{\tiny File:~\jobname .tex,
          printed: \number\day-\number\month-\number\year,
          \thehours.\ifnum\theminutes<10{0}\fi\theminutes}
} \makeatletter\def\thefootnote{\@arabic\c@footnote}\makeatother

\author{J. Chen}
\address{J. Chen, Department of Mathematics,
Hunan Normal University, Changsha, Hunan 410081, People's Republic
of China.} \email{jiaolongchen@sina.com}

\author{ A. Rasila }
\address{A. Rasila, Department of Mathematics,
Hunan Normal University, Changsha, Hunan 410081, People's Republic
of China, and
Department of Mathematics and Systems Analysis, Aalto University, P. O. Box 11100, FI-00076 Aalto,
 Finland.} \email{antti.rasila@iki.fi}

\author{X. Wang${}^{~\mathbf{*}}$}
\address{X. Wang, Department of Mathematics,
Hunan Normal University, Changsha, Hunan 410081, People's Republic
of China.} \email{xtwang@hunnu.edu.cn}


\begin{abstract}
In this paper, we introduce a class of complex-valued polyharmonic mappings, denoted by $HS_{p}(\lambda)$, and its subclass $HS_{p}^{0}(\lambda)$, where $\lambda\in [0,1]$ is a constant. These classes are natural generalizations of a class of mappings studied by Goodman in 1950's. We generalize the main results of Avci and Z{\l}otkiewicz from 1990's to the classes $HS_{p}(\lambda)$ and $HS_{p}^{0}(\lambda)$, showing that the mappings in $HS_{p}(\lambda)$ are univalent and sense preserving. We also prove that the mappings in $HS_{p}^{0}(\lambda)$ are starlike with respect to the origin, and characterize the extremal points of the above classes.
\end{abstract}


\keywords{polyharmonic mapping, starlikeness, convexity, extremal point\\
${}^{\mathbf{*}}$ Corresponding author}

\maketitle

\section{Introduction }\label{csw-sec1}

A complex-valued mapping $F=u+iv$, defined in a domain $D \subset \mathbb{C}$, is called {\it
polyharmonic} (or {\it $p$-harmonic}) if $F$ is $2p$ $(p\geq1)$ times continuously differentiable, and it satisfies the polyharmonic equation $\Delta^{p}F =\Delta(\Delta^{p-1}F)= 0$, where $\Delta^{1}:=\Delta$ is the standard complex Laplacian operator
$$\Delta=4\frac{\partial^{2}}{\partial z\partial \overline{z}}:=
\frac{\partial^{2}}{\partial x^{2}}+\frac{\partial^{2}}{\partial y^{2}}. $$
It is well known (see \cite{shpo1, qiwa}) that for a simply connected domain $D$, a mapping
$F$ is polyharmonic if and only if $F$  has the following representation:
$$F(z)=\sum_{k=1}^{p}|z|^{2(k-1)}G_{k}(z),$$ where $G_{k}$ are
complex-valued harmonic mappings in $D$ for $k\in \{1,\cdots,p\}$. Furthermore, the mappings $G_{k}$ can be expressed as the form
\begin{center}
$G_{k} = h_{k} + \overline{g_{k}}$
\end{center}
for $k\in \{1,\cdots,p\}$, where all $h_{k}$ and  $g_{k}$ are analytic in $D$ (see \cite{cl, du}).

Obviously, for $p=1$ (resp. $p=2$), $F$ is a harmonic (resp. biharmonic) mapping. The biharmonic model arises from numerous problems in science and engineering (cf. \cite{ha, kh, la}). However, investigation of biharmonic mappings in the context of the geometric function theory has been started only recently (see \cite{abab, A, ababkh, CPW0, CPW4,CW}). The reader is referred to \cite{shpo1, qiwa} for discussion on polyharmonic mappings, and \cite{cl, du} for the properties of harmonic mappings.

In \cite{av}, Avci and Z{\l}otkiewicz introduced the class $HS$ of univalent harmonic mappings $F$ with the series expansion:
\be\label{eq0.1}F(z)=h(z)+\overline{g(z)}=
z+\sum_{n=2}^{\infty}a_{n}z^{n}+\sum_{n=1}^{\infty}\overline{b_{n}}\overline{z^{n}}\ee
such that
$$\sum_{n=2}^{\infty}n(|a_{n}|+|b_{n}|)\leq 1-|b_{1}|\;\; (0\leq|b_{1}|<1),$$
and the subclass $HC$ of $HS$, where
$$\sum_{n=2}^{\infty}n^{2}(|a_{n}|+|b_{n}|)\leq 1-|b_{1}|\;\; (0\leq|b_{1}|<1).$$
The corresponding subclasses of $HS$ and $HC$ with $b_{1}=0$ are
denoted by $HS^{0}$ and $HC^{0}$, respectively. These two classes constitute a harmonic
counterpart of classes introduced by Goodman \cite{go}. They are useful in studying questions of so-called $\delta$-neighborhoods (Ruscheweyh \cite{ru}, see also \cite{qiwa}) and in constructing explicit $k$-quasiconformal extensions (Fait et al. \cite{fa}).
In this paper, we define polyharmonic analogues $HS_{p}(\lambda)$ and $HS_{p}^{0}(\lambda)$, where $\lambda \in [0,1]$, to the above classes of mappings. Our aim is to generalize the main results of \cite{av} to the mappings of the classes $HS_{p}(\lambda)$ and $HS_{p}^{0}(\lambda)$.

This paper is organized as follows. In Section \ref{csw-sec2}, we discuss the starlikeness and convexity
of polyharmonic mappings in $HS_{p}^{0}(\lambda)$. Our main result, Theorem \ref{thm2.2}, is a generalization of \cite[Theorem 4]{av}. In Section \ref{csw-sec3}, we find the extremal points of the class $HS_{p}^{0}(\lambda)$. The main result of this section is Theorem \ref{thm2.6}, which is a generalization of \cite[Theorem 6]{av}. Finally, we consider convolutions and existence of neighborhoods. The main results in this section are Theorems \ref{thm2.7} and \ref{thm2.8} which are generalizations of \cite[Theorems $7$ $(i)$ and $8$]{av}, respectively.

\section{Preliminaries}

For $r>0$, write $\mathbb{D}_{r} =\{z:|z|<r\}$, and let
$\mathbb{D}:=\mathbb{D}_{1}$, i.e., the unit disk. We use $H_{p}$ to
denote the set of all polyharmonic mappings $F$ in $\mathbb{D}$ with
a series expansion of the following form: \be\label{eq1.1}
F(z)=\sum_{k=1}^{p}|z|^{2(k-1)}\big(h_{k}(z)+\overline{g_{k}(z)}\big)
=\sum_{k=1}^{p}|z|^{2(k-1)}\sum_{n=1}^{\infty}(a_{n,k}z^{n}+\overline{b_{n,k}}\overline{z^{n}})
\ee with $a_{1,1}=1$ and $|b_{1,1}|<1$. Let $H^{0}_{p}$ denote the
subclass of $H_{p}$ for $b_{1,1}=0$ and $a_{1,k}=b_{1,k}=0$ for
$k\in \{2,\cdots,p\}$.

In \cite{qiwa}, J. Qiao and X. Wang introduced the class $HS_{p}$ of polyharmonic mappings $F$ with the form
\eqref{eq1.1} satisfying the conditions
\be\label{eq0001.0}
\begin{cases}
\displaystyle\sum_{k=1}^{p}\sum_{n=2}^{\infty}(2(k-1)+n)(|a_{n,k}|+|b_{n,k}|)\leq 1-|b_{1,1}|-\sum_{k=2}^{p}(2k-1)(|a_{1,k}|+|b_{1,k}|),\\
\displaystyle0\leq|b_{1,1}|+\sum_{k=2}^{p}(|a_{1,k}|+|b_{1,k}|)<1,
\end{cases}\ee
and the subclass $HC_{p}$ of $HS_{p}$, where
\be\label{eq001.0}\begin{cases}
\displaystyle\sum_{k=1}^{p}\sum_{n=2}^{\infty}(2(k-1)+n^{2})(|a_{n,k}|+|b_{n,k}|)\leq 1-|b_{1,1}|-\sum_{k=2}^{p}(2k-1)(|a_{1,k}|+|b_{1,k}|),\\
\displaystyle0\leq|b_{1,1}|+\sum_{k=2}^{p}(|a_{1,k}|+|b_{1,k}|)<1.
\end{cases}\ee
The classess of all mappings $F$ in $H_{p}^{0}$ which are of the form \eqref{eq1.1}, and subject to the conditions \eqref{eq0001.0}, \eqref{eq001.0}, are denoted by $HS^{0}_{p}$, $HC^{0}_{p}$, respectively.

Now we introduce a new class of polyharmonic mappings, denoted by $HS_{p}(\lambda)$, as follows: A mapping $F\in H_{p}$ with the form \eqref{eq1.1} is said to be in $HS_{p}(\lambda)$ if
\be\label{eq1.0}
\begin{cases}
\displaystyle\sum_{k=1}^{p}\sum_{n=2}^{\infty}\big(2(k-1)+n(\lambda n+1-\lambda) \big)(|a_{n,k}|+|b_{n,k}|)\leq2-\sum_{k=1}^{p}(2k-1)(|a_{1,k}|+|b_{1,k}|),\\
\displaystyle1\leq\sum_{k=1}^{p}(2k-1)(|a_{1,k}|+|b_{1,k}|)<2,
\end{cases}
\ee
where $\lambda \in [0,1].$ We denote by $HS^{0}_{p}(\lambda)$  the class consisting of all mappings $F$ in $H_{p}^{0}$, with the form \eqref{eq1.1}, and subject to the condition \eqref{eq1.0}.
Obviously, if $\lambda=0$ or $\lambda=1$, then the class $HS_{p}(\lambda)$ reduces to $HS_{p}$ or $HC_{p}$, respectively. Similarly, if $p-1=\lambda=0$ or $p=\lambda=1$, then $HS_{p}(\lambda)$ reduces to $HS$ or $HC$.

If $$F(z)=\sum_{k=1}^{p}|z|^{2(k-1)}\sum_{n=1}^{\infty}\left(a_{n,k}z^{n}+
\overline{b_{n,k}}\overline{z^{n}}\right)$$ and
$$G(z)=\sum_{k=1}^{p}|z|^{2(k-1)}\sum_{n=1}^{\infty}\left(A_{n,k}z^{n}+
\overline{B_{n,k}}\overline{z^{n}}\right),$$ then the
{\it convolution} $F\ast G$ of $F$ and $G$ is defined to be the mapping
$$(F\ast G)(z)=\sum_{k=1}^{p}|z|^{2(k-1)}\sum_{n=1}^{\infty}\left(a_{n,k}A_{n,k}z^{n}+
\overline{b_{n,k}}\overline{B_{n,k}}\overline{z^{n}}\right),$$
while the {\it integral convolution} is defined by
 $$(F\diamond G)(z)=\sum_{k=1}^{p}|z|^{2(k-1)}\sum_{n=1}^{\infty}
 \left(\frac{a_{n,k}A_{n,k}}{n}z^{n}+
\frac{\overline{b_{n,k}}\overline{B_{n,k}}}{n}\overline{z^{n}}\right).$$
See \cite{PD-1} for similar operators defined on the class of analytic functions.

Following the notation of J. Qiao and X. Wang \cite{qiwa}, we denote the $\delta$-neighborhood of $F$ the set by
\begin{align*}
N_{\delta}(F(z))=&\left\{G(z):\;\sum_{k=1}^{p}\sum_{n=2}^{\infty}
\big(2(k-1)+n\big)(|a_{n,k}-A_{n,k}|+|b_{n,k}-B_{n,k}|)\right.\\
&\left.+\sum_{k=2}^{p}(2k-1)(|a_{1,k}-A_{1,k}|
+|b_{1,k}-B_{1,k}|)+|b_{1,1}-B_{1,1}|\leq\delta\right\},
\end{align*}
where $G(z)=\sum_{k=1}^{p}|z|^{2(k-1)}\sum_{n=1}^{\infty}\left(A_{n,k}z^{n}+
\overline{B_{n,k}}\overline{z^{n}}\right)$ and $A_{1,1}=1$ (see also Ruscheweyh \cite{ru}).

\section{Starlikeness and convexity}\label{csw-sec2}

 We say that a univalent polyharmonic mapping $F$ with $F(0) = 0$ is
{\it starlike} with respect to the origin if the curve $F(re^{i\theta})$ is starlike with respect to the
origin for each $0 < r < 1$.

\begin{prop}{\rm \bf (\cite{po})}\label{pro1}
If $F$ is univalent, $F(0) = 0$ and $\frac{d}{d \theta}\big(\arg F(re^{i\theta})\big) >0$ for
$z = re^{i\theta}\not= 0$, then $F$ is starlike with respect to the origin.
\end{prop}

 A univalent polyharmonic mapping $F$ with $F(0) = 0$
and $\frac{d}{d \theta}F(re^{i\theta}) \not=0$ whenever $0<r<1$,
is said to be {\it convex} if the curve $F(re^{i\theta})$ is convex for each $0 < r < 1$.

\begin{prop}{\rm \bf (\cite{po})} \label{pro2} If $F$ is univalent,
$F(0) = 0$ and $\frac{\partial}{\partial \theta} \left[ \arg \left(\frac{\partial}{\partial \theta}F(re^{i\theta})\right) \right]  >0$ for
$z = re^{i\theta}\not= 0$, then $F$ is convex.
\end{prop}

Let $X$ be a
topological vector space over the field of complex
numbers, and let $D$ be a set of $X$. A point $x \in D$ is
called an {\it extremal point} of $D$ if
it has no representation of the form $x=ty+(1-t)z$ $(0<t<1)$ as a proper convex
combination of two distinct points $y$ and $z$ in $D$.

Now we are ready to prove results concerning the geometric
properties of mappings in $HS_{p}^{0}(\lambda)$.

\begin{Thm}{\rm \bf \cite[Theorems $3.1$, $3.2$ and $3.3$]{qiwa}}\label{k1} Suppose that $F\in HS_{p}$. Then $F$ is
univalent and sense preserving in $\mathbb{D}$. In particular, each
member of $HS_{p}^{0}$ (or $HC_{p}^{0}$) maps $\mathbb{D}$ onto a
domain starlike w.r.t. the origin, and a convex domain,
respectively.
\end{Thm}

\bthm\label{thm2.2} Each mapping in $HS^{0}_{p}(\lambda)$ maps the disk $\mathbb{D}_{r}$, where
$r\leq\max\{\frac{1}{2},\lambda\}$, onto a convex domain. \ethm

\bpf Let $F\in HS^{0}_{p}(\lambda)$, and let $r\in (0,1)$ be fixed. Then
 $r^{-1}F(rz)\in HS^{0}_{p}(\lambda)$ by \eqref{eq1.0}, and we have
\begin{align*}
&\sum_{k=1}^{p}\sum_{n=2}^{\infty}\big(2(k-1)+n^{2} \big)(|a_{n,k}|+|b_{n,k}|)r^{2k+n-3}\\
\leq&\sum_{k=1}^{p}\sum_{n=2}^{\infty}\big(2(k-1) +n(\lambda n
+1-\lambda)\big)(|a_{n,k}|+|b_{n,k}|)\leq1
\end{align*}
provided that $$\big(2(k-1)+n^{2} \big)r^{2k+n-3}\leq2(k-1)
+n(\lambda n +1-\lambda)$$ for $k\in\{1,\cdots,p\}$, $n\geq2$
and $0\leq\lambda\leq1$, which is true if $r\leq\max\{\frac{1}{2},\lambda\}$.
Then the result follows from Theorem A.
 \epf

Follows immediately from Theorem A, we get the following.
\bcor\label{thm2.1} Let $F\in HS_{p}(\lambda)$. Then $F$ is a
univalent, sense preserving polyharmonic mapping. In particular, if
$F\in HS_{p}^{0}(\lambda)$, then $F$ maps $\mathbb{D}$ onto a domain
starlike w.r.t. the origin. \ecor

\begin{figure}
\centering
\includegraphics[width=2.5in]{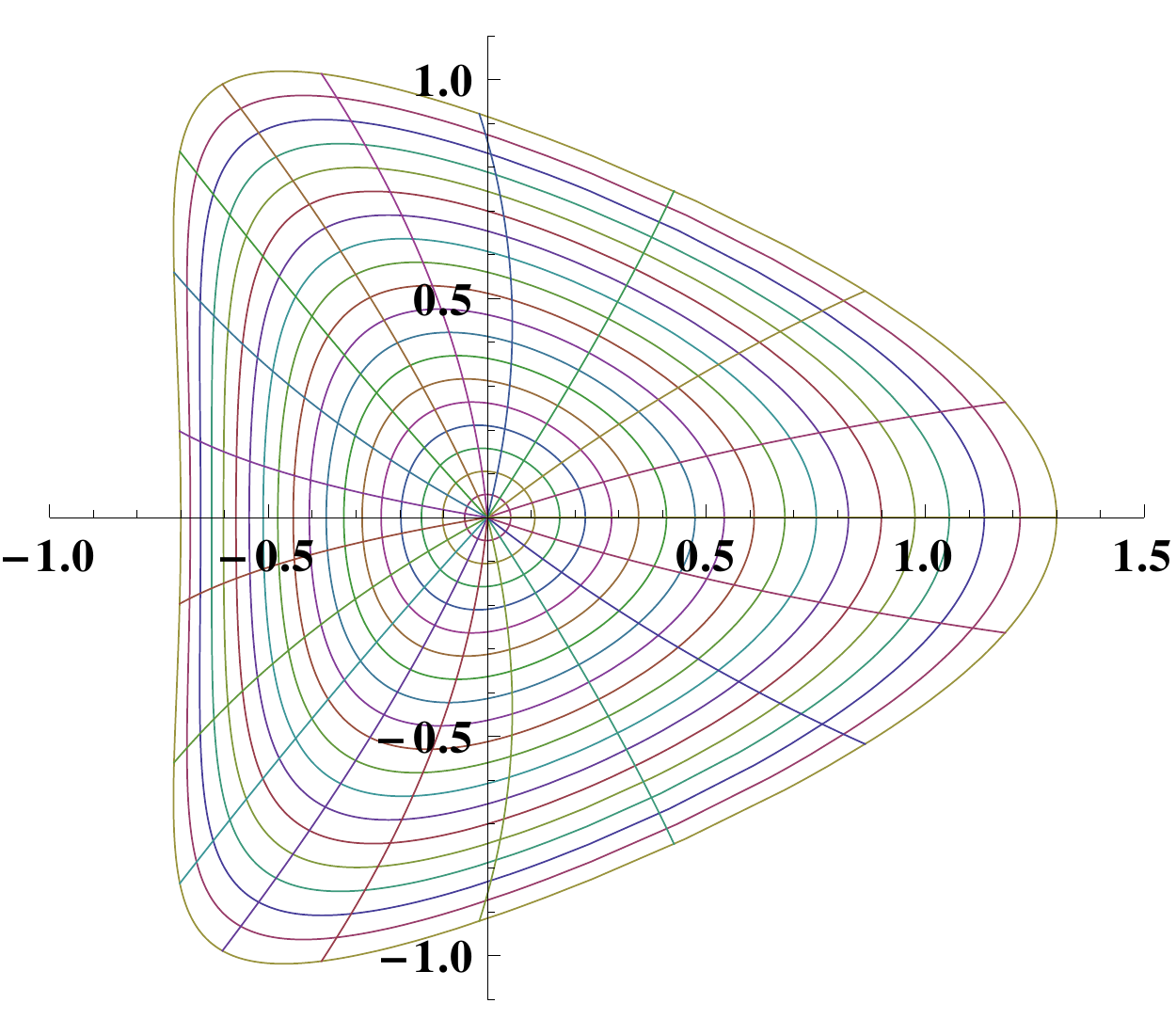}
\qquad
\includegraphics[width=2.2in]{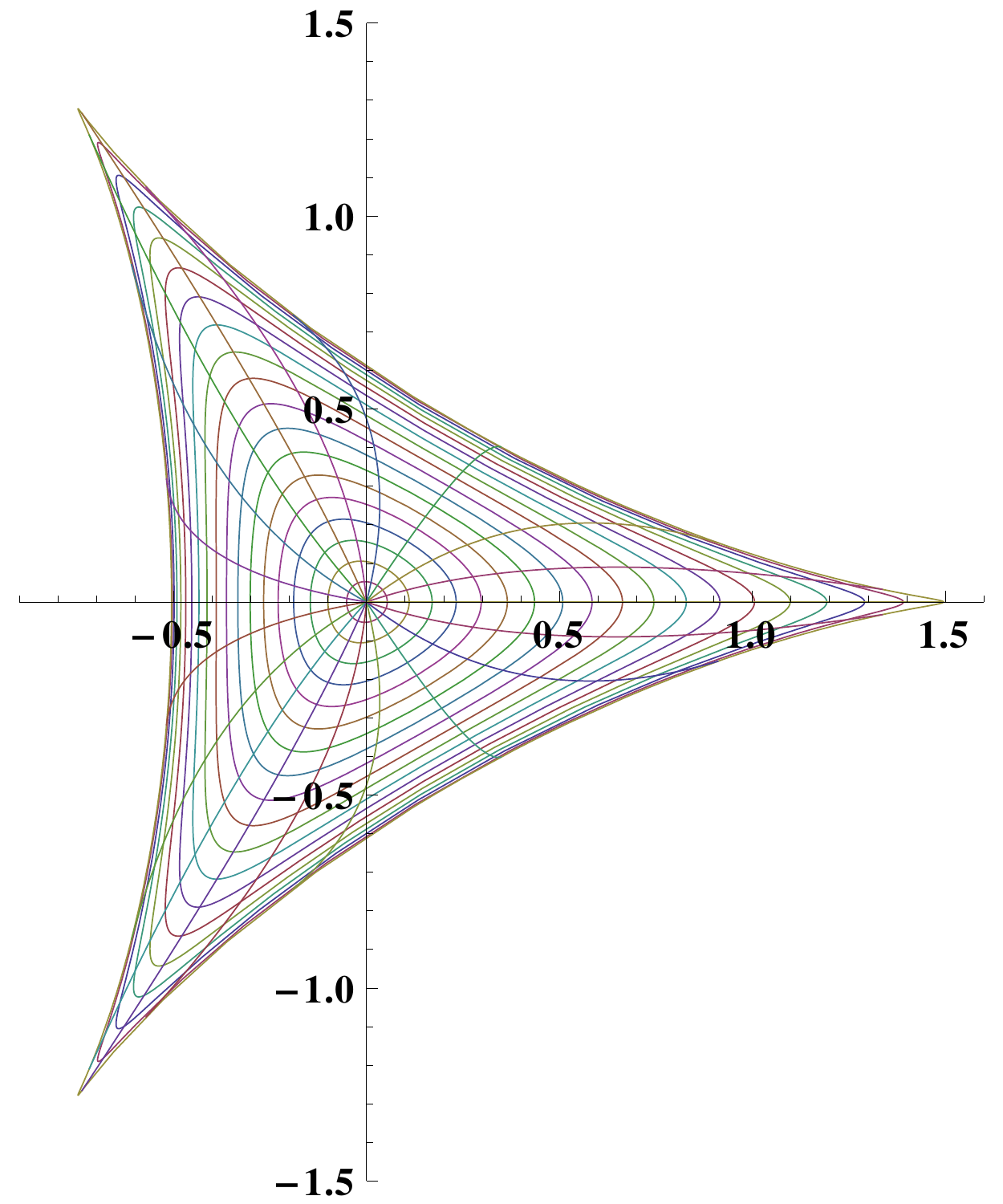}
\caption{The images of $\mathbb{D}$ under the mappings $F_{1}(z)=z+\frac{1}{10}z^2+\frac{1}{5}\overline{z^{2}}$ (left) and $F_{2}(z)=z+\frac{1}{101}z^2+\frac{49}{101}\overline{z^{2}}$ (right). }\label{f2}
\end{figure}

\begin{examp} Let $F_{1}(z)=z+\frac{1}{10}z^{2}+\frac{1}{5}\overline{z^{2}}$.
Then $F_{1}\in HS_{1}^{0}(\frac{2}{3})$ is a univalent, sense
preserving polyharmonic mapping. In particular, $F_{1}$ maps $\mathbb{D}$
onto a domain starlike w.r.t. the origin, and it maps the disk $\mathbb{D}_{r}$,
where $r\leq\frac{2}{3}$, onto a convex domain. See Figure \ref{f2}.
\end{examp}This example shows that the class $HS_{p}^{0}(\lambda)$ of polyharmonic mappings is more general than the class $HS^{0}$
which is studied in \cite{av} even in the case of harmonic mappings (i.e. $p$=1).

\begin{examp} Let $F_{2}(z)=z+\frac{1}{101}z^2+\frac{49}{101}\overline{z^{2}}$.
Then $F_{2}\in HS_{1}^{0}(\frac{1}{100})$ is a univalent, sense
preserving polyharmonic mapping. In particular, $F_{2}$ maps $\mathbb{D}$
onto a domain starlike w.r.t. the origin, and it maps the disk $\mathbb{D}_{r}$, where
$r\leq\frac{1}{2}$, onto a convex domain. See Figure \ref{f2}.
\end{examp}

\section{Extremal points}\label{csw-sec3}
First, we determine the distortion bounds for mappings in $HS_{p}(\lambda)$.
\blem\label{thm2.3} Suppose that $F\in HS_{p}(\lambda)$.
Then the following statements hold:

$(1)$ For $0\leq\lambda\leq\frac{1}{2}$,
$$ (1-|b_{1,1}|)|z|-\frac{1-|b_{1,1}|}{2(1+\lambda)}|z|^{2}\leq
|F(z)|\leq(1+|b_{1,1}|)|z|+\frac{1-|b_{1,1}|}{2(1+\lambda)}|z|^{2}.$$
Equalities are obtained by the mappings
$$F(z)=z+|b_{1,1}|e^{i\mu}\overline{z}+\frac{1-|b_{1,1}|}{2(1+\lambda)}e^{i\nu}z^{2},$$
for properly chosen real $\mu$ and $\nu$;

$(2)$ For $\frac{1}{2}<\lambda\leq1$,
$$|F(z)\leq(1+|b_{1,1}|)|z|+\frac{1-|b_{1,1}|-3(|a_{1,2}|
+|b_{1,2}|)}{2(1+\lambda)}|z|^{2}+(|a_{1,2}|+|b_{1,2}|)|z|^{3}$$
and $$|F(z)|\geq(1-|b_{1,1}|)|z|-\frac{1-|b_{1,1}|-3(|a_{1,2}|
+|b_{1,2}|)}{2(1+\lambda)}|z|^{2}-(|a_{1,2}|+|b_{1,2}|)|z|^{3}.$$
Equalities are obtained by the mappings
 $$F(z)=z+|b_{1,1}|e^{i\eta}\overline{z}+
 \frac{1-|b_{1,1}|-3(|a_{1,2}|+|b_{1,2}|)}{2(1+\lambda)}e^{i\varphi}z^{2}
 +(|a_{1,2}|+|b_{1,2}|)e^{i\psi}z|z|^{2},$$
 for properly chosen real $\eta$, $\varphi$ and $\psi$.
\elem
\bpf Let $F\in HS_{p}(\lambda)$, where $\lambda\in[0,1]$. By \eqref{eq1.1}, we have
$$
|F(z)|\leq(1+|b_{1,1}|)|z|+\left(\sum_{k=1}^{p}\sum_{n=2}^{\infty}(|a_{n,k}|+|b_{n,k}|)
+\sum_{k=2}^{p}(|a_{1,k}|+|b_{1,k}|)
\right)|z|^{2}.
$$
For $0\leq\lambda\leq\frac{1}{2}$, we have \be\label{eq1.19}
2(1+\lambda)\leq 2k-1, \ee where $k\in \{2,\cdots,p\}$,
 and
\be\label{eq1.20}
 2(1+\lambda)\leq 2(k-1)+n(\lambda n+1-\lambda),
 \ee where $k\in \{1,\cdots,p\}$ and $n\geq2$.
Then \eqref{eq1.19}, \eqref{eq1.20} and \eqref{eq1.0} give
\begin{align*}
&\sum_{k=1}^{p}\sum_{n=2}^{\infty}(|a_{n,k}|+|b_{n,k}|)
+\sum_{k=2}^{p}(|a_{1,k}|+|b_{1,k}|)\\
\leq&\frac{1}{2(1+\lambda)}\Big(1-|b_{1,1}|-\sum_{k=1}^{p}\sum_{n=2}^{\infty}
\big(2(k-1)+n(\lambda n+1-\lambda)-2(1+\lambda)  \big)(|a_{n,k}|+|b_{n,k}|)\\
&-\sum_{k=2}^{p}\big((2k-1)-2(1+\lambda) \big)(|a_{1,k}|+|b_{1,k}|)\Big),\\
\end{align*}
so $$(1-|b_{1,1}|)|z|-\frac{1-|b_{1,1}|}{2(1+\lambda)}|z|^{2}
\leq|F(z)|\leq(1+|b_{1,1}|)|z|+\frac{1-|b_{1,1}|}{2(1+\lambda)}|z|^{2}.$$
By \eqref{eq1.1}, we obtain

\beq \nonumber |F(z)| &\leq&
(1+|b_{1,1}|)|z|+\left(\sum_{k=1}^{p}\sum_{n=2}^{\infty}(|a_{n,k}|+|b_{n,k}|)
+\sum_{k=3}^{p}(|a_{1,k}|+|b_{1,k}|)\right)|z|^{2}\\ \nonumber
&+&(|a_{1,2}|+|b_{1,2}|)|z|^{3}.\eeq For $\frac{1}{2}<\lambda\leq1$,
we have \be\label{eq1.21} 2(1+\lambda)\leq 2k-1, \ee  where $k\in
\{3,\cdots,p\}$,
 and
\be\label{eq1.22}
 2(1+\lambda)\leq 2(k-1)+n(\lambda n+1-\lambda),
 \ee where $k\in \{1,\cdots,p\},\;n\geq 2$.
Then \eqref{eq1.21}, \eqref{eq1.22} and \eqref{eq1.0} imply
\begin{align*}
&\sum_{k=1}^{p}\sum_{n=2}^{\infty}(|a_{n,k}|+|b_{n,k}|)
+\sum_{k=3}^{p}(|a_{1,k}|+|b_{1,k}|)\\
\leq&\frac{1}{2(1+\lambda)}\Big(1-|b_{1,1}|-\sum_{k=1}^{p}\sum_{n=2}^{\infty}
\big(2(k-1)+n(\lambda n+1-\lambda)-2(1+\lambda)  \big)(|a_{n,k}|+|b_{n,k}|)\\
&-\sum_{k=3}^{p}\big(2k-1-2(1+\lambda) \big)
(|a_{1,k}|+|b_{1,k}|)-3(|a_{1,2}|+|b_{1,2}|)\Big).\\
\end{align*}
Then $$|F(z)|\geq(1-|b_{1,1}|)|z|-\frac{1-|b_{1,1}|
-3(|a_{1,2}|+|b_{1,2}|)}{2(1+\lambda)}|z|^{2}-(|a_{1,2}|+|b_{1,2}|)|z|^{3}$$
and
$$|F(z)|\leq(1+|b_{1,1}|)|z|+\frac{1-|b_{1,1}|
-3(|a_{1,2}|+|b_{1,2}|)}{2(1+\lambda)}|z|^{2}+(|a_{1,2}|+|b_{1,2}|)|z|^{3}.$$
The proof of this lemma is complete.\epf

\br\label{rem1} Suppose that $F\in HS_{p}(\lambda)$ is of the form
$$F(z)=\sum_{k=1}^{p}|z|^{2(k-1)}G_{k}(z)=\sum_{k=1}^{p}|z|^{2(k-1)}\sum_{n=1}^{\infty}\left(a_{n,k}z^{n}+
\overline{b_{n,k}}\overline{z^{n}}\right).$$
Then for each $k\in \{1,\cdots,p\}$,
$$|G_{k}(z)|\leq  (|a_{1,k}|+|b_{1,k}|)|z|+\frac{1-|b_{1,1}|}{2(1+\lambda)}|z|^{2}.$$
\er

\blem\label{thm2.4} The family $HS_{p}(\lambda)$ is closed under
convex combinations.
\elem

\bpf Suppose $F_{i}\in HS_{p}(\lambda)$ and  $t_i\in [0,1]$
with $\sum_{i=1}^{\infty}t_i=1$. Let
$$
F_{i}(z)=\sum_{k=1}^{p}|z|^{2(k-1)}\sum_{n=1}^{\infty}
(a^{(i)}_{n,k}z^{n}+\overline{b^{(i)}_{n,k}}\overline{z^{n}}).$$ By
Lemma \ref{thm2.3}, there exists a constant $M$ such that
$|F_{i}(z)|\leq M$ for all $i=1,\cdots,p$. It follows that
$\sum_{i=1}^{\infty}t_iF_{i}(z)$ is absolutely and uniformly
convergent, and by Remark \ref{rem1}, the mapping
$\sum_{i=1}^{\infty}t_i F_{i}(z) $ is polyharmonic. Since
$\sum_{i=1}^{\infty}t_iF_{i}(z)$ is absolutely and uniformly
convergent, we have
\begin{align*}
\sum_{i=1}^{\infty}t_iF_{i}(z)&=\sum_{i=1}^{\infty}t_i\sum_{k=1}^{p}
|z|^{2(k-1)}\left(\sum_{n=1}^{\infty}a^{(i)}_{n,k}z^n+\sum_{n=1}^{\infty}\overline{b^{(i)}_{n,k}}\overline{z^n}\right)\\
&=\sum_{k=1}^p|z|^{2(k-1)}\left(\sum_{n=1}^{\infty}\sum_{i=1}^{\infty}t_ia^{(i)}_{n,k}z^n
+\sum_{n=1}^{\infty}\sum_{i=1}^{\infty}t_i \overline{b^{(i)}_{n,k}}\overline{z^n}\right).\\
\end{align*}
By \eqref{eq1.0}, we get
\be\label{eq1.30}
\sum_{k=1}^{p}\sum_{n=1}^{\infty}\big(2(k-1)+n(\lambda n+1-\lambda) \big)
\left(\left|\sum_{i=1}^{\infty}t_ia^{(i)}_{n,k}\right|+\left|\sum_{i=1}^{\infty}t_i b^{(i)}_{n,k}\right|\right)\ee
$$\leq\sum_{i=1}^{\infty}t_i\left ( \sum_{k=1}^{p}\sum_{n=1}^{\infty} \big(2(k-1)+n(\lambda  n+1-\lambda)\big)
(|a^{(i)}_{n,k}|+|b^{(i)}_{n,k}|)  \right)\leq2.$$ It follows from
$$1\leq\sum_{k=1}^{p}(2k-1)
\left(\left|\sum_{i=1}^{\infty}t_ia^{(i)}_{1,k}\right|+
\left|\sum_{i=1}^{\infty}t_ib^{(i)}_{1,k}\right|\right)<2$$ and
\eqref{eq1.30} that $\sum_{i=1}^{\infty}t_iF_{i}\in
HS_{p}(\lambda)$. \epf

From Lemma \ref{thm2.3}, we see that the class $HS_{p}(\lambda)$ is
uniformly bounded, and hence normal. Lemma \ref{thm2.4} implies
that $ HS^{0}_{p}(\lambda)$ is also compact and convex. Then there exists
a non-empty set of extremal points in $HS_{p}^{0}(\lambda)$.

\bthm\label{thm2.6}
The extremal points of $HS_{p}^{0}(\lambda)$ are the mappings
of the following form:
$$
F_k(z)=z+|z|^{2(k-1)}a_{n,k}z^n\;\text{or}\;
F^*_k(z)=z+|z|^{2(k-1)}\overline{b_{m,k}}\overline{z^m},
$$
where
$$|a_{n,k}|=\frac{1}{2(k-1)+n(\lambda n+1-\lambda)},\;\;\text{for} \;n\geq2,\;k\in\{1,\cdots,p\},$$
and
$$|b_{m,k}|=\frac{1}{2(k-1)+m(\lambda m+1-\lambda)},\;\;\text{for} \;m\geq2,\;k\in\{1,\cdots,p\}.$$
 \ethm

\bpf Assume that $F$ is an extremal point of $HS_{p}^{0}(\lambda)$, of the form $(\ref{eq1.1})$.
Suppose that the coefficients of $F$ satisfy the following:
$$
\sum_{k=1}^{p}\sum_{n=2}^{\infty}\big(2(k-1)+n(\lambda n+1-\lambda) \big)(|a_{n,k}|+|b_{n,k}|)<1.
$$
If all coefficients $a_{n,k}$ $(n\geq2)$ and $b_{n,k}$ $(n\geq2)$
are equal to $0$, we let
$$
F_1(z)=z+\frac{1}{2(1+\lambda)}z^2\;\text{and} \;
F_2(z)=z-\frac{1}{2(1+\lambda)}z^2.
$$
Then $F_1$ and $F_2$ are in $ HS_{p}^{0}(\lambda)$ and
$F=\frac{1}{2}(F_1+F_2)$. This is a contradiction, showing that
there is a coefficient, say $a_{n_{0},k_{0}}$ or $b_{n_{0},k_{0}}$,
of $F$ which is nonzero. Without loss of generality, we may further
assume that $a_{n_{0},k_{0}}\neq0$.

For $\gamma>0$ small enough, choosing $x\in \mathbb{C}$ with $|x|=1$ properly and
replacing $a_{n_{0},k_{0}}$ by $a_{n_{0},k_{0}}-\gamma x$ and
$a_{n_{0},k_{0}}+\gamma x$, respectively, we obtain two mappings
$F_3$ and $F_4$ such that both $F_3$ and $F_4$ are in $HS_{p}^{0}(\lambda)$. Obviously, $F=\frac{1}{2}(F_3+F_4)$.
Hence the coefficients of $F$ must satisfy the following equality:
$$
\sum_{k=1}^{p}\sum_{n=2}^{\infty}\big(2(k-1)+n(\lambda n+1-\lambda) \big)(|a_{n,k}|+|b_{n,k}|)=1.
$$

Suppose that there exists at least two coefficients, say,
$a_{q_1,k_1}$ and $b_{q_2,k_2}$ or $a_{q_1,k_1}$ and
$a_{q_2,k_2}$ or $b_{q_1,k_1}$ and $b_{q_2,k_2}$, which
are not equal to 0, where $q_1$, $q_2\geq2$. Without loss of generality, we assume the first
case. Choosing $\gamma>0$ small enough and  $x\in \mathbb{C}$, $y\in \mathbb{C}$ with
$|x|=|y|=1$ properly, leaving all coefficients of $F$ but
$a_{q_1,k_1}$ and $b_{q_2,k_2}$ unchanged and replacing
$a_{q_1,k_1}, b_{q_2,k_2}$ by
$$
a_{q_1,k_1}+\frac{\gamma x}{2(k_{1}-1)+q_{1}
(\lambda q_{1}+1-\lambda)}\;\;\text{and}\;\;b_{q_2,k_2}
-\frac{\gamma y}{2(k_{2}-1)+q_{2}(\lambda q_{2}+1-\lambda)},
$$
or
$$
a_{q_1,k_1}-\frac{\gamma x}{2(k_{1}-1)+q_{1}
(\lambda q_{1}+1-\lambda)}\;\;\text{and}\;\;b_{q_2,k_2}
+\frac{\gamma y}{2(k_{2}-1)+q_{2}(\lambda q_{2}+1-\lambda)},
$$
respectively, we obtain two mappings $F_5$ and $F_6$ such that $F_5$
and $F_6$ are in $HS_{p}^{0}(\lambda)$. Obviously,
$F=\frac{1}{2}(F_5+F_6)$. This shows that any extremal
point $F\in HS_{p}^{0}(\lambda)$ must have the form
$F_k(z)=z+|z|^{2(k-1)}a_{n,k}z^n$
or $F^*_k(z)=z+|z|^{2(k-1)}\overline{b_{m,k}}\overline{z^m}$,
where
$$|a_{n,k}|=\frac{1}{2(k-1)+n(\lambda n+1-\lambda)},\;\;\text{for} \;n\geq2,\;k\in\{1,\cdots,p\},$$
and
$$|b_{m,k}|=\frac{1}{2(k-1)+m(\lambda m+1-\lambda)},\;\;\text{for} \;m\geq2,\;k\in\{1,\cdots,p\}.$$

Now we are ready to prove that for any $F\in HS_{p}^{0}(\lambda)$
with the above form must be an extremal point of
$HS_{p}^{0}(\lambda)$. It suffices to prove the case of
$F_k$, since the proof for the case of $F^*_k$ is similar.

Suppose there exist two mappings $F_7$ and $F_8\in
HS_{p}^{0}(\lambda)$ such that $F_k=tF_7+(1-t)F_8\;(0<t<1).$
For $q=7,8$, let
$$
F_{q}(z)=\sum_{k=1}^p|z|^{2(k-1)}\sum_{n=1}^{\infty}(a^{(q)}_{n,k}z^n
+\overline{b^{(q)}_{n,k}}\overline{z^n}).
$$
Then \be\label{eq2.2}
|ta_{n,k}^{(7)}+(1-t)a_{n,k}^{(8)}|=|a_{n,k}|=\frac{1}{2(k-1)+n(\lambda n+1-\lambda)}.
\ee

Since all coefficients of $F_q\;(q=7,8)$ satisfy, for $n\geq2$ and $k\in\{1,\cdots, p\}$,
$$
|a_{n,k}^{(q)}|\leq
\frac{1}{2(k-1)+n(\lambda n+1-\lambda)},\quad |b_{n,k}^{(q)}|\leq
\frac{1}{2(k-1)+n(\lambda n+1-\lambda)},
$$
(\ref{eq2.2}) implies $a_{n,k}^{(7)}=a_{n,k}^{(8)}$, and all other
coefficients of $F_7$ and $F_8$ are equal to $0$. Thus
$F_k=F_7=F_8$, which shows that $F_k$ is an extremal point of
$HS_{p}^{0}(\lambda).$
 \epf

\section{Convolutions and neighborhoods}\label{csw-sec4}

Let $C_{H}^{0}$ denote the class of harmonic univalent, convex
mappings $F$ of the form \eqref{eq0.1} with $b_{1}=0$.
It is known \cite{cl} that the below sharp inequalities hold:
$$2|a_{n}|\leq n+1,\;\;2|b_{n}|\leq n-1.$$

It follows from \cite[Theorems $5.14$]{cl} that if $H$ and $G$ are in $C_{H}^{0}$,
then $H\ast G$ (or $H\diamond G$) is sometime not convex,
but it may be univalent or even convex if one of the mappings $H$ and $F$ satisfies some additional conditions.
In this section, we consider convolutions of harmonic mappings
$F\in HS_{1}^{0}(\lambda)$ and $H\in C^{0}_{H}$.

\bthm\label{thm2.7} Suppose that $H(z)=z+\sum_{n=2}^{\infty}(A_{n}z^{n}
+\overline{B_{n}}\overline{z^{n}})\in C_{H}^{0}$
and $F\in HS_{1}^{0}(\lambda)$. Then for $\frac{1}{2}\leq \lambda\leq 1$, the convolution $F\ast H$ is
univalent and starlike, and the integral convolution $F\diamond H$ is convex.
\ethm

\bpf If $F(z)=z+\sum_{n=2}^{\infty}(a_{n}z^{n}+\overline{b_{n}}\overline{z^{n}})
\in HS_{1}^{0}(\lambda)$, then for $F\ast H$, we obtain
\begin{align*}
\sum_{n=2}^{\infty}n(|a_{n}A_{n}|+|b_{n}B_{n}|)\leq&\sum_{n=2}^{\infty}n
\left(\frac{n+1}{2}|a_{n}|+\frac{n-1}{2}|b_{n}|\right)\\
\leq&\sum_{n=2}^{\infty}n(\lambda n+1-\lambda)(|a_{n}|+|b_{n}|)\leq1.
\end{align*}
Hence $(F\ast H)\in HS^{0}$. The transformations
$$\int_{0}^{1}\frac{F(z)\ast H(tz)}{t} dt=(F\diamond H)(z)$$
now show that $(F\diamond H)\in HC^{0}$. By Theorem A, the result follows.
\epf

\br\label{rem22} The proof of the Theorem \ref{thm2.7} does not generalize to polyharmonic mappings, when $p\ge 2$. For example, let $p=2$, and write $$H(z)=z+\sum_{k=1}^{2}|z|^{2(k-1)}\sum_{n=2}^{\infty}(A_{n,k}z^{n}
+ \overline{B_{n,k}}\overline{z^{n}})$$ and
$$F(z)=z+\sum_{k=1}^{2}|z|^{2(k-1)}\sum_{n=2}^{\infty}(a_{n,k}z^{n}
+ \overline{b_{n,k}}\overline{z^{n}}).$$
Suppose that $|A_{n,k}|\leq\frac{n+1}{2}$, $|B_{n,k}|\leq\frac{n-1}{2}$ and $F\in HS_{2}^{0}(\lambda)$.
Then for $\lambda=1$, the convolution $F\ast H$ is
univalent and starlike but it is not clear if this is true for $\frac{1}{2}\le \lambda <1$. However, the integral convolution $F\diamond H$ is convex for $\frac{1}{2}\leq \lambda \leq1$.
\er

\begin{figure}
\centering
\includegraphics[width=2.2in]{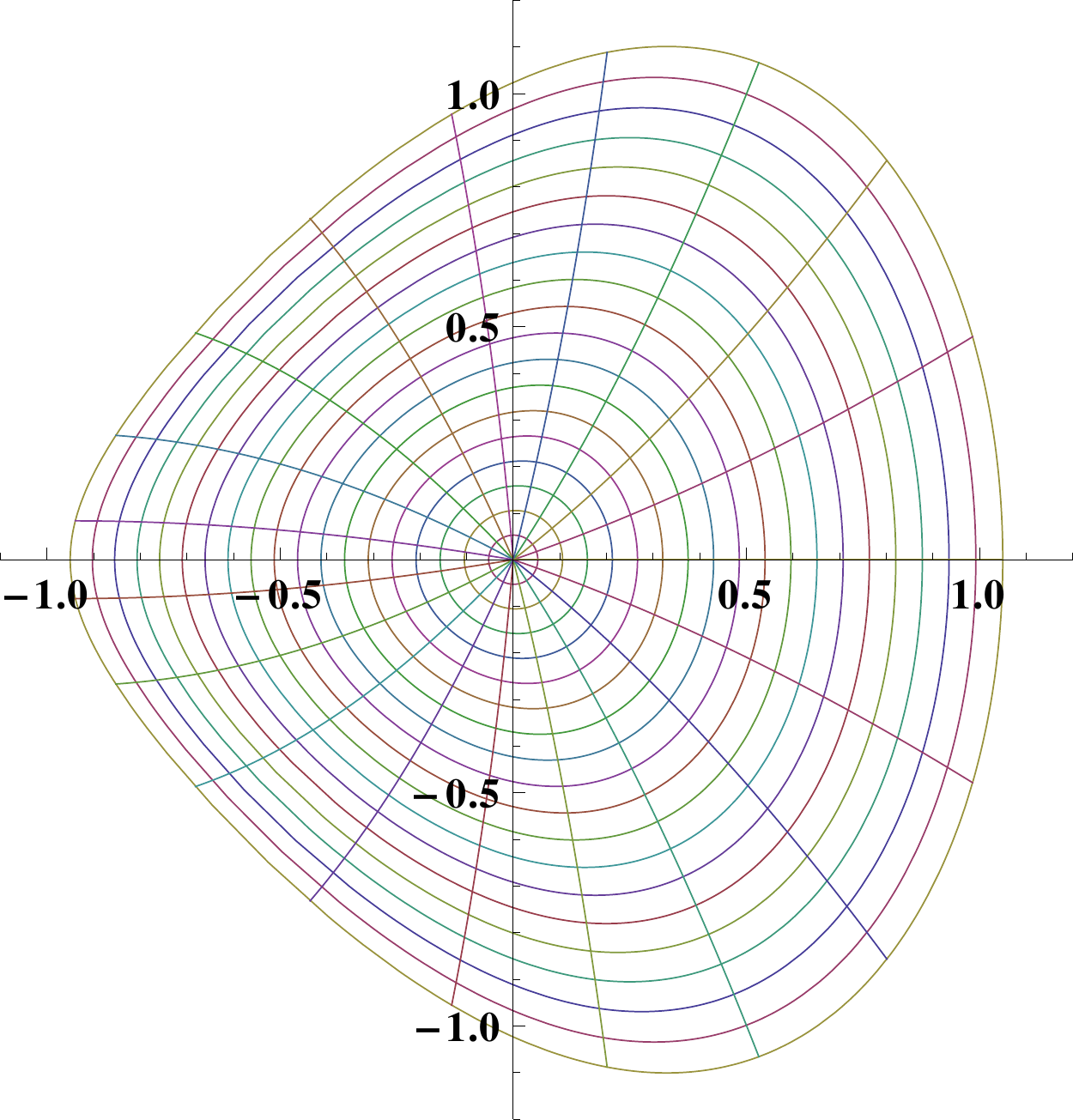}
\qquad
\includegraphics[width=2.3in]{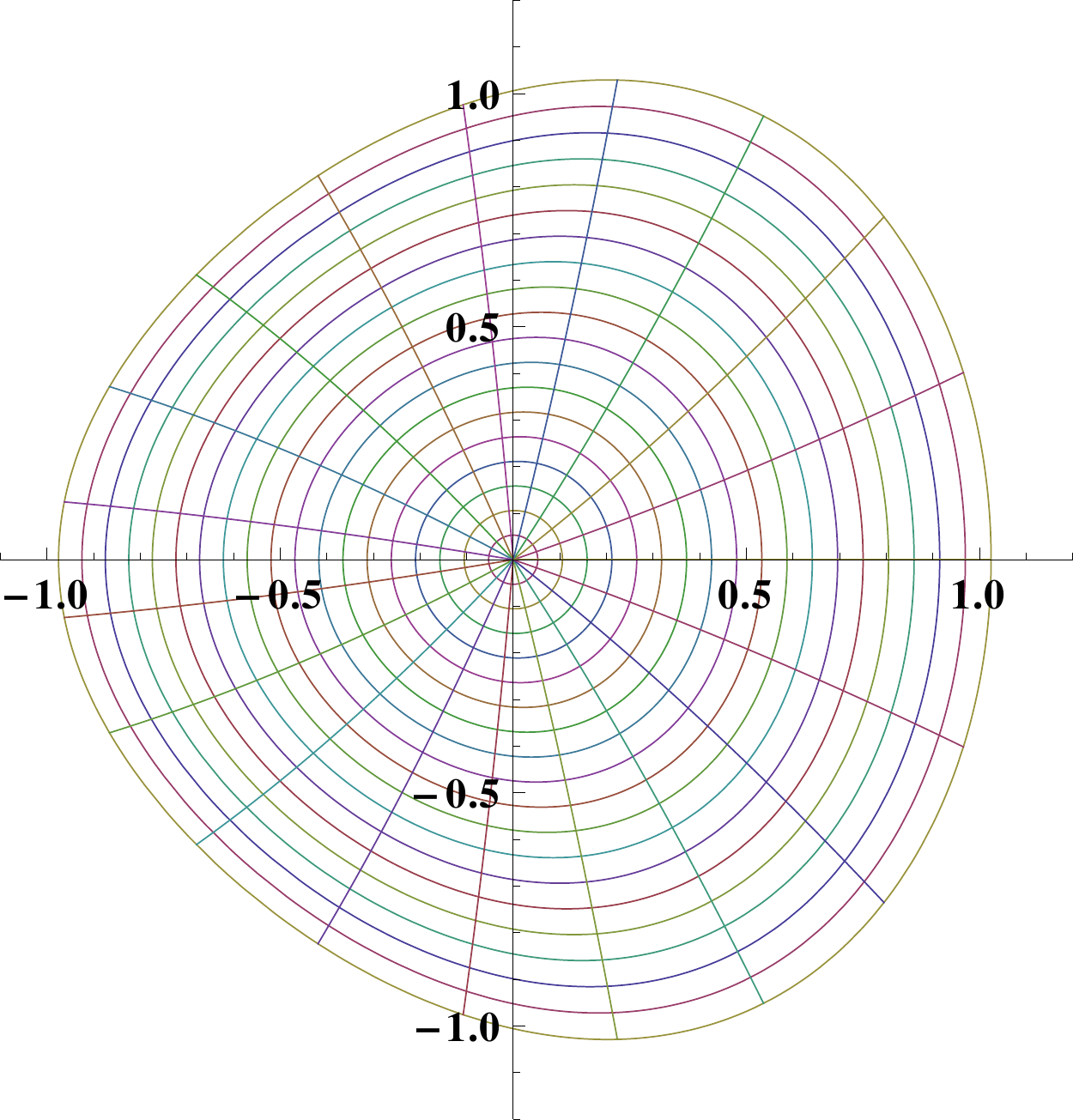}
\caption{The images of $\mathbb{D}$ under the mappings $(F\ast H)(z)=z+\frac{3}{20}z^{2}-\frac{1}{10}\overline{z^{2}}$ (left) and $(F\diamond H)(z)=z+\frac{3}{40}z^{2}-\frac{1}{20}\overline{z^{2}}$ (right).}\label{f5}
\end{figure}

\begin{examp} Let $H(z)$=$\rm{Re}$$\big\{\frac{z}{1-z} \big\}+i$$\rm{Im}$$\big\{\frac{z}{(1-z)^{2}}
\big\}\in$$C_{H}^{0}$. Then $H(z)$ maps $\mathbb{D}$ onto the
half-plane $\rm{Re}$$\{w\}>\frac{1}{2}$, and let
$F(z)=z+\frac{1}{10}z^{2}+\frac{1}{5}\overline{z^{2}}\in
HS^{0}_{1}(\frac{2}{3})$. Then the convolution $F\ast H$ is
univalent and starlike, and the integral convolution $F\diamond H$
is convex (see Figure \ref{f5}).
\end{examp}

Finally, we are going to prove the existence of neighborhoods for mappings in the class $HS_{p}(\lambda)$.

\bthm\label{thm2.8} Assume that $\lambda \in (0,1]$ and $F\in HS_{p}(\lambda)$. If $$\delta\leq\frac{\lambda}{p+\lambda}
\left(2-\sum_{k=1}^{p}(2k-1)(|a_{1,k}|+|b_{1,k}|) \right),$$ then $N_{\delta}(F)\subset HS_{p}$.
\ethm

\bpf Let $H(z)=\sum_{k=1}^{p}|z|^{2(k-1)}\sum_{n=1}^{\infty}
(A_{n,k}z^{n}+\overline{B_{n,k}}\overline{z}^{n})\in N_{\delta}(F)$. Then
\begin{align*}
&\sum_{k=1}^{p}\sum_{n=2}^{\infty}\big(2(k-1)+n\big)(|A_{n,k}|+|B_{n,k}|)
+\sum_{k=2}^{p}(2k-1)(|A_{1,k}|+|B_{1,k}|)+|B_{1,1}|\\
\leq&\sum_{k=1}^{p}\sum_{n=2}^{\infty}\big(2(k-1)+n\big)(|A_{n,k}-a_{n,k}|+|B_{n,k}-b_{n,k}|)\\
&+\sum_{k=2}^{p}(2k-1)(|A_{1,k}-a_{1,k}|+|B_{1,k}-b_{1,k}|)+|B_{1,1}-b_{1,1}|\\
&+\sum_{k=1}^{p}\sum_{n=2}^{\infty}\big(2(k-1)+n\big)(|a_{n,k}|+|b_{n,k}|)
+\sum_{k=2}^{p}(2k-1)(|a_{1,k}|+|b_{1,k}|)+|b_{1,1}|\\
\leq&\delta+\sum_{k=1}^{p}\sum_{n=2}^{\infty}\big(2(k-1)+n\big)
(|a_{n,k}|+|b_{n,k}|)+\sum_{k=2}^{p}(2k-1)(|a_{1,k}|+|b_{1,k}|)+|b_{1,1}|\\
\leq&\delta+\frac{p}{p+\lambda}\sum_{k=1}^{p}\sum_{n=2}^{\infty}\big(2(k-1)+n(\lambda n+1-\lambda)\big)(|a_{n,k}|+|b_{n,k}|)\\&+\sum_{k=2}^{p}(2k-1)(|a_{1,k}|+|b_{1,k}|)+|b_{1,1}|\\
\leq&\delta+\frac{p-\lambda}{p+\lambda}+\frac{\lambda}{p+\lambda}\sum_{k=1}^{p}(2k-1)(|a_{1,k}|+|b_{1,k}|)
<1.
\end{align*}Hence, $H\in HS_{p}$.
\epf

\end{document}